\documentclass{article}
\usepackage[utf8]{inputenc}

\title{Acyclic Comprehension is equal to Stratified Comprehension}
\author{Zuhair Al-Johar and M. Randall Holmes}
\date{June 8 2010, Updated: January 22 2011\footnote{This document is a type-set copy of the original document with simple modifications, copied at October 7, 2020.}}

\begin{document}

\maketitle

\section{Introduction}
A new criterion of comprehension is defined, initially termed by myself as ``connected” and finally 
as ``Acyclic" by Mr. Randall Holmes. Acyclic comprehension simply asserts that for any acyclic formula $\phi$ the set $\{x : \phi\} $ exists. I first presented this criterion semi-formally to Mr. Randall Holmes, who further made the first rigorous definition of it, a definition that I finally simplified to the one presented here. Later Mr. Holmes made another presentation of the definition which is also mentioned here. He pointed to me that acyclic comprehension is implied by stratification, and posed the question as to whether it is equivalent to full stratification or strictly weaker. He and initially I myself thought that it was strictly weaker; Mr. Randall Holmes actually conjectured that it is very weak. Surprisingly it turned to be equivalent to full stratification as I proved here.\bigskip

We are indebted to Mr. Nathan Bowler for singling out an error with the prior proof of inclusion

\section {Definition of Acyclic Comprehension }

We say that a variable $x$ is connected to a variable $y$ in the formula $\phi$ iff any of the following formulae appear in $\phi : x \in y; \ y \in x; \ x = y; \ y = x.$ \bigskip

We refer to a function $s$ from $\{1, ..., n\}$ to variables in $\phi$ as a chain of length $n$ in $\phi$ iff for each appropriate index $i: s_i$ is connected to $s_{i+1}$, and for each appropriate index $j: s_j , s_{j+2}$ are two different occurrences in $\phi$. A chain from $x$ to $y$ is defined as a chain $s$ of length $n > 1$ with $s_1 = x$ and $s_n = y$.\bigskip

A formula $\phi$ is said to be acyclic iff for each variable $x$ in $\phi$ there is no chain from $x$ to $x$.\bigskip

\noindent
\textbf{\emph{Mr. Randall Holmes definition of acyclic formulae:}}

With any formula $\phi$ associate a non-directed graph $G_\phi$ whose vertices are the variables in $\phi$ and which contain an edge from $x$ to $y$ for each atomic formula 
$x \in y; \ y \in x; \ x = y; \ y = x$ which $occurs$ as a subformula of $\phi$.\bigskip

$\phi$ is said to be acyclic iff $G_\phi$ is acyclic.\bigskip

\textbf{\emph{Acyclic Comprehension:}} For $n = 0, 1, 2, ...$; if $\phi$ is acyclic formula in first order logic with identity and membership, in which $y$ is free, and in which $x$ does not occur, then: $\forall w_1...w_n \ \exists x \ \forall y \ (y \in x \Leftrightarrow \phi)$

\section{Stratified Comprehension implies acyclic Comprehension}

Suppose that $\phi$ is acyclic. From each component $C$ of the graph $G_\phi$ choose a variable $x_c$. Assign type 0 to each $x_c$. For each variable $y$ in the component $C$, there is a unique path from $y$ to $x_c$, in the multigraph $G_\phi$. Assign types for each variable occurring in this path consistent with the requirements for stratification of atomic formulae, and a type is assigned to $y$. The uniqueness of the path ensures that this is actually an effective definition of a stratification. (R.Holmes).

\section {Acyclic Comprehension implies stratified Comprehension}

The proof relies upon proving that all axioms of the finite axiomatization of stratified
comprehension present in the second author's online book [1] as adapted to Wiener’s ordered pairs [2], are proved by acyclic comprehension and Extensionality (weak or full). The argument that this is 
a finite axiomatization of stratified comprehension is the same one behind adaptation to Kuratowski’s ordered pairs explained in [1]. The bi-conditional is taken to be primitive.

Unless otherwise noted, strong Extensionality is assumed throughout the following proof:\bigskip

\textbf {Universal Set:} $V \equiv_{def}  \{y|\exists s: s = y\}$ \bigskip

\textbf{Complements:} $a^c \equiv_{def}  \{y|y \not \in  a\}$\bigskip

\textbf{Boolean union:} $a \cup  b \equiv_{def}  \{y|y \in a \lor y \in b\}$\bigskip

\textbf {Set union:} $\bigcup(a) \equiv_{def}  \{y|\exists z(y \in z \in a)\}$\bigskip

\textbf {Singletons:} $\{a\} \equiv_{def}  \{y|y = a\}$\bigskip

\textbf {Wiener’s ordered pairs:} $(a, b) \equiv_{def}  \{\{\{a\}, 0\}, \{\{b\}\}\}$\bigskip

$(a, b) \equiv_{def}  \{y|\forall w(w \in y \Leftrightarrow \forall k(k \in w \Leftrightarrow k = a) \lor w = 0)\lor$

$\ \ \ \ \ \ \ \ \ \ \ \ \ \ \ \ \ \ \ \forall u(u \in y \Leftrightarrow \forall n(n \in u \Leftrightarrow n = b))\}$\bigskip

Note: if weak Extensionality is assumed, then $0$ is to be replaced by $V$.

(R. Holmes).\bigskip

$(, )(x) \equiv_{def}  \exists a, b(x = (a, b))$\bigskip

\textbf {Cartesian products:} $A \times B \equiv_{def}  \{x|\exists a, b (x = (a, b) \land  a \in A \land  b \in B)\}$\bigskip

\textbf {Domains:} $Dom(R) \equiv_{def}  \{y|\exists m \in R \exists p, q(y \in p \in q \in m \land  0 \in q)\}$\bigskip

\textbf {Diagonal:} $ [=] \equiv_{def}  \{(i, i)|i \in V \} \equiv_{def}  \{y|(, )(y) \land  \exists !i\exists p(i \in p \in y \land  i \neq 0)\}$\bigskip

\textbf {Converses:}
$Conv(R) \equiv_{def}  \{(b, a)|(a, b) \in R\} \equiv_{def}  \{y|(, )(y) \land$  

$\exists m \in R \forall q, l, p, k((p \in q \in m \land  k \in l \in y \land  (0 \not \in q \Leftrightarrow 0 \in l) \land  p \neq 0 \land  k \neq 0) \Rightarrow $

$ p = k)\}$\bigskip

\textbf {Singleton images:}
$SI(R) \equiv_{def}  \{(\{a\}, \{b\})|(a, b) \in R\} \equiv_{def}  \{y|(, )(y) \land  $

$\exists  m \in R  \forall q, l, p, k((p \in q \in
m \land  k \in l \in y \land  (0 \in q \Leftrightarrow 0 \in l) \land  p \neq 0 \land  k \neq 0) \Rightarrow$

$p \in k)\}$\bigskip

\textbf {Relative Products:}
$R|S \equiv_{def}  \{(a, b)|\exists c ((a, c) \in R \land  (c, b) \in S)\}$\bigskip

$\exists !x:\phi$  is  to be  written  as   $``\exists z\forall x(\phi \Leftrightarrow x = z)"$, were  $z$  do  not  occur  in  $\phi$ \bigskip

$C(R \times S) \equiv_{def}  \{((a, c), (c, b))|(a, c) \in R \land  (c, b) \in S\}$ \bigskip

$\equiv_{def}  \{y|y \in R \times S \land  \exists !i\exists g, h, k, l(i \in g \in h \in k \in l \in y \land  i \neq 0 \land  (0 \not \in g \Leftrightarrow$

$ 0 \in l))\}$\bigskip

$R|S \equiv_{def}  \{y|(, )(y) \land  \exists m \in C(R \times S)\forall f, g, h, k, l,w, u((f \in g \in h \in k \in $

$l \in m \land  w \in u \in y \land  (0 \in u \Leftrightarrow 0 \in l \Leftrightarrow 0 \in g) \land  f \neq 0 \land  w \neq 0) \Rightarrow  f = w)\}$.\bigskip

\textbf {Projections:}

$\pi_1 \equiv_{def}  \{((x, y), \{\{\{x\}\}\})|x \in V \land  y \in V \}$ \smallskip

$\pi_2 \equiv_{def}  \{((x, y), \{\{\{y\}\}\})|x \in V \land  y \in V \}$\bigskip

$V_1 \equiv_{def}  \{(a, b)|(, )(a) \land  \exists u(b = \{\{\{u\}\}\})\}$\smallskip

$V_2 \equiv_{def}  \{y|y \in V_1 \land  \exists z\forall w(\exists o, p(w \in o \in p \in y) \Rightarrow  \exists r, s(z \in r \in s \in w))\}$\smallskip

$V! \equiv_{def}  \{y|y \in V_1 \land  \exists !t\exists b, n,m, o(t \in b \in n \in m \in o \in y \land  t \neq 0)\}$\bigskip

$\pi_2 \equiv_{def}  \{y|y \in V_2 \in \land  \exists !j\exists r, s, p, q(j \in r \in s \in p \in q \in y \land  0 \not \in r)\}$\smallskip

$\pi_1 \equiv_{def}  (V_2 \backslash \pi_2) \cup V!$\bigskip

\textbf {Inclusion:}
$[\subseteq] \equiv_{def}  \{(z,x)| z \subseteq x\}\equiv _{def}  \bigcup(\{t| \exists x\forall y (y \in t \Leftrightarrow \forall w (\exists m\in w \land $

$\exists k\forall u (u \in w \Leftrightarrow
((,)(u) \land \forall i (\exists s, r (i \in s \in r \in u) \Rightarrow  i \subseteq k) \land  \forall j (\exists p, q (j \in p \in$ 

$  q \in u \land  0 \not \in q )\Rightarrow  j = x))) \Rightarrow 
y \in w))\})$ \bigskip

\noindent
Note: if weak Extensionality is assumed, then technical modification of some proofs is required in order to allow for the replacement of $0$ by $V$ in Wiener’s ordered pairs.

\section{References}

[1]Holmes, M. R. [1998] 
Elementary set theory with a universal set. volume 10 of the Cahiers du Centre de logique, Academia, Louvain-la-Neuve (Belgium), 241 
pages, ISBN 2-87209-488-1. Official On-line version at:

\noindent
http://math.boisestate.edu/~holmes/holmes/head.pdf \smallskip

\noindent
[2] Wiener, Norbert, “A simplification of the logic of relations”, in From Frege to
Gödel, a sourcebook in mathematical logic, 1879-1931 (van Heijenoort, editor),
Harvard University Press 1967, pp. 224-227. \bigskip

\end{document}